\documentclass[11pt,reqno]{amsart}
\usepackage{amsmath,amssymb,amsthm}
\usepackage{graphicx}
\usepackage[all]{xypic}
\usepackage[utf8,nocaptions]{vietnam}
\usepackage{calrsfs}
\input xypic
\def\thmsection{section}
\def\thmchangesection{changesection}

\def\thmchangechapter{changechapter}
\def\thmchange{change}
\def\thmplain{plain}
\ifx\theoremnumberstyle\thmplain
  \theoremstyle{break-italic}
  \newtheorem{satz}{Satz}
\else
  \ifx\theoremnumberstyle\thmsection
    \theoremstyle{break-italic}
    \newtheorem{satz}{Satz}[section]
  \else
    \ifx\theoremnumberstyle\thmchange
      \swapnumbers
      \theoremstyle{break-italic}
      \newtheorem{satz}{Satz}
    \else
      \ifx\theoremnumberstyle\thmchangesection
         \swapnumbers
         \theoremstyle{break-italic}
         \newtheorem{satz}{Satz}[section]
      \else
        \ifx\theoremnumberstyle\thmchangechapter
           \swapnumbers
           \theoremstyle{break-italic}
           \newtheorem{satz}{Satz}[chapter]
        \else
          \ifx\theoremnumberstyle\thmsection
             \theoremstyle{break-italic}
             \newtheorem{satz}{Satz}[section]
          \else
            \theoremstyle{break-italic}
            \newtheorem{satz}{Satz}[section]
          \fi
        \fi
      \fi
    \fi
  \fi
\fi

\theoremstyle{break-italic}
\newtheorem{theorem}[satz]{Theorem}
\newtheorem{lemma}[satz]{Lemma}

\newtheorem{corollary}[satz]{Corollary}
\newtheorem{Proposition}[satz]{Proposition}

\newtheorem*{conjecture*}{Conjecture}

\theoremstyle{break-roman}
\newtheorem{definition}[satz]{Definition}

\newtheorem{example}[satz]{Example}

\newtheorem{remark}[satz]{Remark}

\newtheorem{conjecture}[satz]{Conjecture}

\theoremstyle{standard}

\newtheorem*{claim}{Claim}

\theoremstyle{varthm-roman}
\newtheorem*{varthm-roman}{}
\theoremstyle{varthm-italic}
\newtheorem*{varthm-italic}{}
\theoremstyle{varthm-roman-break}
\newtheorem*{varthm-roman-break}{}
\theoremstyle{varthm-italic-break}
\newtheorem*{varthm-italic-break}{}
\theoremstyle{varthm-roman-no-punctuation}
\newtheorem{varthm-roman-no-punctuation-numbered}[satz]{}
\theoremstyle{varthm-italic-no-punctuation}
\newtheorem{varthm-italic-no-punctuation-numbered}[satz]{}

\newenvironment{varthm-roman-numbered}[1]{
  \begin{varthm-roman-no-punctuation-numbered}
    \mbox{\rm\textbf{#1}}
  }{\end{varthm-roman-no-punctuation-numbered}}
\newenvironment{varthm-italic-numbered}[1]{
  \begin{varthm-italic-no-punctuation-numbered}
    \mbox{\rm\textbf{#1}}
  }{\end{varthm-italic-no-punctuation-numbered}}
\newenvironment{varthm-roman-break-numbered}[1]{
  \begin{varthm-roman-no-punctuation-numbered}
    \mbox{\rm\textbf{#1}\newline}
  }{\end{varthm-roman-no-punctuation-numbered}}
\newenvironment{varthm-italic-break-numbered}[1]{
  \begin{varthm-italic-no-punctuation-numbered}
    \mbox{\rm\textbf{#1}}\newline
  }{\end{varthm-italic-no-punctuation-numbered}}
\numberwithin{equation}{section}

\def\ex{\begin{example}
  }
  \def\eex{\end{example}}
\def\thr{\begin{theorem}}
\def\ethr{\end{theorem}}
\def\pro{\begin{Proposition}}
\def\epro{\end{Proposition}}
\def\coro{\begin{corollary}}
\def\ecoro{\end{corollary}}
\def\df{\begin{definition}}
\def\edf{\end{definition}}
\def\lm{\begin{lemma}}
\def\elm{\end{lemma}}
\def\pf{\begin{proof}}
\def\epf{\end{proof}}
\def\problem{\begin{problem}}
\def\eproblem{\end{problem}}

\def\ite{\begin{itemize}}
\def\hite{\end{itemize}}
\def\rem{\begin{remark}}
\def\erem{\end{remark}}

\def\cla{\begin{claim}}
\def\ecla{\end{claim}}
\def\conj{\begin{conjecture}}
\def\econj{\end{conjecture}}

\begin{document}

\title[A note on optimization with Morse polynomials]{A note on optimization with Morse polynomials}
\author{L\^{e} C\^{o}ng-Tr\`{i}nh }
\address{L\^{e} C\^{o}ng-Tr\`{i}nh \\
Department of Mathematics, Quy Nhon University\\
170 An Duong Vuong, Quy Nhon, Binh Dinh}
\email{lecongtrinh@qnu.edu.vn}

\date{\today}


\begin{abstract} In this paper we prove that the gradient ideal of a Morse polynomial is radical. This gives  a \textit{generic class} of polynomials whose gradient ideals are radical. As a consequence we reclaim a previous result  that the unconstrained polynomial  optimization problem for Morse polynomials  has a finite convergence. 
\end{abstract}
\maketitle


\section{Introduction}

Let us consider the \textit{unconstrained   polynomial  optimization problem}
\begin{equation}\label{pt-global}
f^*=\inf_{x\in \mathbb R^n} f(x),
\end{equation}
 where $f\in \mathbb R[x_1,\cdots,x_n]=:\mathbb R[X]$ is a real polynomial in $n$ variables $x_1,\cdots, x_n$ whose degree is $d$. 
 
 In principle polynomial optimization problems can  be solved algorithmically by symbolic computation using algorithms from real algebraic geometry \cite{BPR03, S08}.  Moreover, the optimization  problem (\ref{pt-global}) is NP-hard even when $d=4$ (cf. \cite{Ne00}).
 
 A milestone in minimizing multivariate polynomials was given by J.-B. Lasserre (2001, \cite{La01}). He has introduced a \textit{sum of squares \footnotemark \footnotetext{A \textit{sum of squares} in $\mathbb R[X]$ is a sum of the form $\displaystyle \sum_{i=1}^k f_i^2$, where $k\in \mathbb N$ and $f_i\in \mathbb R[X]$ for every $i=1,\cdots, k$.}  (SOS) relaxation} for this problem by computing a lower bound for $f^*$: 
\begin{align}\label{pt-sos}
 f^*_{sos}&: = \mbox{ maximize }  \gamma \\
 & \mbox{ subject to } f(x) - \gamma \mbox{ is a SOS in } \mathbb R[X].\nonumber
 \end{align}
 
Note that the SOS relaxation (\ref{pt-sos}) can be solved by semidefinite programming (SDP). Moreover,  it should be also noted that 
\begin{align*}
 f^* & = \mbox{ maximize }  \gamma  \\
 & \mbox{ subject to } f(x) - \gamma \geq 0 \mbox{ for all } x \in \mathbb R^n.
 \end{align*}
Therefore  $f^*_{sos}\leq f^*$, and the equality occurs if and only if $f-f^*$ is a SOS. In the case where $f^*_{sos}$ is strictly less than $f^*$, Lasserre  proposed finding a sequence of lower bounds for $f(x)$ in some large ball $\{x\in \mathbb R^n | \quad  \|x\|\leq R\}$. This sequence converges to $f^*$ when the degrees of the polynomials introduced in the algorithm tend to infinity. However, it may not converge in finitely many steps, and moreover the degrees of the required auxiliary polynomials can be very large. 

J. Nie, J. Demmel and B. Sturmfels (2006, \cite{NDS06}) have given a method to compute the global infimum $f^*$ which terminates in finitely many steps in terms of  the  gradient ideal $I_{grad}(f)=\left<\frac{\partial f}{\partial x_1},\cdots,\frac{\partial f}{\partial x_n} \right>$.  The finite convergence of their SOS relaxation depends  mainly on the SOS representation of $f-f^{*}$ modulo the gradient ideal $I_{grad}(f)$.  More explicitly, for each natural number $N$, let us consider the  problem
\begin{equation}\label{grad-relax}
f^*_{N, grad}:=\mbox{maximize \quad }   \gamma
\end{equation}
$ \mbox{ subject to } f(x)-\gamma - \sum_{i=1}^m \phi_i(x)\frac{\partial f }{\partial x_i}  \mbox{ is a SOS in } \mathbb R[X] \mbox{  and }   \phi_i(x) \in \mathbb R[X]_{2N-d+1}.$ 

Here $\mathbb R[X]_m$ denotes the $\binom{n+m}{m}$-dimensional vector space of polynomials of degree at most $m$. Nie, Demmel and Sturmfels proved the following 
\thr[{\cite[Theorem 10]{NDS06}}] \label{thrNDS06} If $f$ attains its infimum $f^*$ over $\mathbb R^n$, then the sequence  $\{f^*_{N,grad}\}_N$ converges increasing to $f^*$.  Moreover, if the gradient ideal $I_{grad}(f)$ is radical, then there exists an integer $N$ such that $f^*_{N,grad}=f^*$.
\ethr
It follows from Theorem \ref{thrNDS06} that the radicality of the gradient ideal $I_{grad}(f)$ is a sufficient condition for the finite convergence of the sequence $\{f^*_{N,grad}\}$. Moreover, for \textit{almost all} polynomials $f$ in $\mathbb R[X]_d$ the gradient ideal $I_{grad}(f)$ is radical \cite[Proposition 1]{NDS06}, i.e. the ideal $I_{grad}(f)$ is generically radical. 

In this paper we introduce  a \textit{generic class} of polynomials whose gradient ideals are radical - the class of \textit{Morse polynomials}.  Recall that a point $x^0\in \mathbb R^n$ is called  a \textit{critical point} of a polynomial $f\in \mathbb R[X]$ if 
$$ \dfrac{\partial f }{\partial x_1}(x^0)=\cdots =  \dfrac{\partial f }{\partial x_n}(x^0)=0. $$
  A critical point $x^0$ of the  polynomial $f $   is said to be \textit{ non-degenerate  }  if the Hessian matrix\footnotemark \footnotetext{For a polynomial $f\in \mathbb R[x_1,\cdots,x_n]$, its \textit{Hessian matrix} at a point $x^0\in \mathbb R^n$ is defined by 
   $$  \nabla^2 f(x^0):= \left[\begin{matrix}
\frac{\partial^2 f}{\partial x_1^2}(x^0) & \cdots & \frac{\partial^2 f}{\partial x_1 \partial x_n}(x^0)\\
\vdots &  & \vdots \\
\frac{\partial^2 f}{\partial x_n\partial x_1}(x^0) & \cdots & \frac{\partial^2 f}{\partial x_n^2}(x^0)
\end{matrix}\right].$$ } $\nabla^{2}f(x^0)$ of $f$ at $x^0$ is invertible.  A polynomial $f\in \mathbb R[X]$ is called a  \textit{Morse polynomial } if all of its critical points are non-degenerate with distinct critical values. It is well-known from differential topology and singularity theory that almost all polynomial  functions on $\mathbb R^n$ are Morse (cf. \cite{BH04}). Moreover, it was shown recently by H\`{a} and Phạm \cite[Theorem 5.1]{HP17} that the set of Morse polynomials in $\mathbb{R}[X]$ whose  \textit{Newton polyhedra at infinity}\footnotemark \footnotetext{ A subset $\mathcal{N}$ of the positive cone $\mathbb R_+^n$ is called a \textit{Newton polyhedron at infinity} if it is the convex hull in $\mathbb R^n$ of the set $A\cup \{0\}$, where  $A \subseteq \mathbb N^n$ is a finite set.} $\mathcal{N}(f)$ contained in a   convenient \footnotemark \footnotetext{A polyhedron $\mathcal{N}\subseteq \mathbb R_+^n$ is called \textit{convenient} if it intersects each coordinate axis at a point different from the origin. 
 For a polynomial  $f(x)=\sum_{\alpha}u_\alpha x^\alpha:=\sum_{\alpha}u_\alpha x_1^{\alpha_1}\cdots x_n^{\alpha_n}\in \mathbb R[X]$, its  \textit{Newton polyhedron at infinity } $\mathcal{N}(f)$ is defined by  the convex hull in $\mathbb R^n$ of the set $\{\alpha | u_\alpha \not = 0\}\cup \{0\}$.  }  Newton polyhedron at infinity $\mathcal{N}$ is open and dense in $\mathcal{N}$. Hence the class of Morse polynomials in $\mathbb R[X]$ with a given convenient Newton polyhedron at infinity is large enough to consider for optimization. 
 
 The main result in this paper is the following 
\thr \label{thr-main} Let $f\in \mathbb R[X]$ be a Morse polynomial. Then its gradient ideal $I_{grad}(f)$ is radical.
\ethr    

As a direct consequence of this theorem and Theorem \ref{thrNDS06}, we get the following result which was obtained by the author in \cite{Le16}.
\coro \label{coro-main} Let $f\in \mathbb R[X]$ be a Morse polynomial. Assume that $f$ attains its infimum $f^*$ over $\mathbb R^n$. Then  there exists an integer $N$ such that $f^*_{N,grad}=f^*$.
\ecoro 

\rem \rm In Corollary \ref{coro-main},  the assumption  that $f$ achieves a minimum on $\mathbb R^n$  cannot be removed. In fact, let us consider the polynomial $f(x,y)=x^2+(xy-1)^2 \in \mathbb R[x,y]$. It is easy to check that $f$ has  a  single critical point at $(0,0)$ which is non-degenerate. Hence $f$ is a Morse polynomial. Moreover, since $\displaystyle\lim_{y\rightarrow \infty}f(\frac{1}{y},y)=0$, we have $f^*=0$. But $f$ does not attain the value $0$ at any points in $\mathbb R^2$. On the other hand, we have 
$$ f(x,y)-1 = \frac{1}{2}x(1-y^2)\frac{\partial f }{ \partial x} +  \frac{1}{2}y(1+y^2)\frac{\partial f }{ \partial y}.$$
Therefore, it follows  from the definition  of $f^*_{N,grad}$ in (\ref{grad-relax}), with $\phi_1(x,y)= \dfrac{1}{2}x(1-y^2)$ and 
$\phi_2(x,y)= \dfrac{1}{2}y(1+y^2)$,  that    $f^*_{N,grad}=1>0=f^*$ for every $N\geq 3$. This means that the sequence  $\{f^*_{N,grad}\}_N$ never  converges   to $f^*$ in this case. 
\erem

\section{Proof of Theorem \ref{thr-main}}

In this section we give a proof for Theorem \ref{thr-main}, applying some known results in singularity theory.  Let  $f\in \mathbb R[X]$ be a Morse polynomial.  Let us denote by
$$ V_{grad}(f):=\{x\in \mathbb C^n |  \frac{\partial f}{\partial x_1}(x)=\cdots=\frac{\partial f}{\partial x_n}(x)=0\} $$ 
the set of \textit{complex} critical points of $f$. Since $f$ is a polynomial, the set $V_{grad}(f)$ is actually an algebraic subset of $\mathbb C^n$ defined by the gradient ideal $I_{grad}(f)$. 
\lm \label{lm-finite} If $f$ is a Morse polynomial, then   $V_{grad}(f)$ is finite.
\elm
\pf It is well-known that each complex critical point of the Morse polynomial $f$ is isolated (cf. \cite[Corollary 2.3]{Mil63}). Therefore the set $V_{grad}(f)$ is discrete. On the other hand,  the set $V_{grad}(f)$ is algebraic in $\mathbb C^n$, hence it has only finitely many connected components (cf. \cite[p.408]{Lo91}). Note that, each critical point itself is a connected component of the algebraic set $V_{grad}(f)$. This implies that $V_{grad}(f)$ has only finitely many elements. 
\epf 
Assume $V_{grad}(f)=\{p^1,\cdots, p^N\}$, and let $p\in \mathbb C^n$ denote  one of the points $p^i$. If $p=(p_1,\cdots,p_n)\in \mathbb C^n$, denote by $\mathbb C\{x_1-p_1,\cdots,x_n-p_n\}$ the ring of convergent power series at the point $p$. The number 
$$\mu(f,p):=\dim_\mathbb C \mathbb C\{x_1-p_1,\cdots,x_n-p_n\}/I_{grad}(f) $$
is called the \textit{Milnor number  of $f$ at $p$}. 

\lm \label{mu-finite} $\mu(f,p)$ is finite.

\elm 
\pf  The result follows from the Hilbert's Nullstellensatz, since $p$ is an isolated critical point of $f$  (cf. \cite[Lemma I. 2.3]{GLS07}). 
\epf 

Denote
$$\mu(f):=\sum_{i=1}^N \mu(f,p^i)=\dim_\mathbb C  \mathbb C[x_1,\cdots,x_n]/I_{grad}(f), $$
  which is called the \textit{total Milnor number} of $f$.  It follows from Lemma \ref{mu-finite} that  \textit{$\mu(f)$ is finite}.

The following result gives a  relation between the number  $\left| V_{grad}(f)\right|$ of complex critical points of $f$  and the total Milnor number of $f$, which yields also a criterion to check for radicality of the gradient ideal $I_{grad}(f)$.
\lm[{\cite[Theorem 2.10]{CLO98}}] \label{lm-radical}   $\left |  V_{grad}(f) \right| \leq \dim_\mathbb C \mathbb C[x_1,\cdots,x_n]/I_{grad}(f)$, with   equality if and only if the gradient ideal $I_{grad}(f)$ is radical.
\elm

To apply this lemma, we need  to prove the following
\lm \label{lm-one}$\mu(f,p^i)=1$ for every $i=1,\cdots, N$.
\elm 
\pf Let us denote by $p$ one of the points $p^i$, $i=1,\cdots, N$. By
translating the variety accordingly, without loss of generality we may assume that $p$ is the
origin $0\in \mathbb C^n$. Consider the   polynomial map 
$$ (f_1,\cdots,f_n): \mathbb C^n \rightarrow \mathbb C^n, $$
where $f_i:=\frac{\partial f }{ \partial x_i} \in \mathbb C[x_1,\cdots,x_n]$, $ i=1,\cdots,n$. We have $f_i(0)=0$ for all $ i=1,\cdots,n$. Moreover, we will show in the following that the determinant of the Jacobian matrix $J(f_1,\cdots,f_n)$ of this  polynomial map at the origin $0\in \mathbb C^n$  is non-zero. In fact, since $f$ is   Morse, we see that $0$ is a non-degenerate critical point of $f$. It follows that the determinant of the Hessian matrix  $\displaystyle \nabla^2f(0)$ of $f$ at $0$ is non-zero.  Observe that  
\begin{align*} 
\det J(f_1,\cdots,f_n)(0) &= \det \left[\begin{matrix}
\frac{\partial f_1}{\partial x_1}(0) & \cdots & \frac{\partial f_1}{\partial x_n}(0)\\
\vdots &  & \vdots \\
\frac{\partial f_n}{\partial x_1}(0) & \cdots & \frac{\partial f_n}{\partial x_n}(0)
\end{matrix}\right] \\
&= \det\left[\begin{matrix}
\frac{\partial^2 f}{\partial x_1^2}(0) & \cdots & \frac{\partial^2 f}{\partial x_1 \partial x_n}(0)\\
\vdots &  & \vdots \\
\frac{\partial^2 f}{\partial x_n\partial x_1}(0) & \cdots & \frac{\partial^2 f}{\partial x_n^2}(0)
\end{matrix}\right]\\
&=\det \nabla^2 f(0).
\end{align*}
Hence $\det J(f_1,\cdots,f_n)(0) \not = 0$. It follows from the \textit{implicit function theorem} for convergent power series  (cf. \cite[Theorem I.1.18]{GLS07}) the following isomorphism of $\mathbb C$-vector spaces:
$$ \mathbb C\{x_1,\cdots,x_n\}/\left<f_1,\cdots,f_n \right> \cong \mathbb C.  $$
This implies that $\mu(f,0)=\dim_\mathbb C \mathbb C =1$.
\epf 
\pf[ \textbf{Proof of Theorem \ref{thr-main}}]  Assume that $f$ is a Morse polynomial. Then it follows from Lemma \ref{lm-finite}  that   $V_{grad}(f)$ is a  finite set  whose elements are $p^1,\cdots, p^N$. By Lemma \ref{lm-one}, at each complex critical point $p^i$ of $f$  we have $\mu(f,p^i)=1$. Hence 
$$\mu(f)=\sum_{i=1}^N \mu(f,p^i)=N=\left|V_{grad}(f)\right|.$$
The theorem now follows from Lemma \ref{lm-radical}.
\epf

\hspace{-0.6cm} \textbf{Acknowledgements~~}  The author would like to express his gratitude to Professor H\`{a} Huy Vui  for his  valuable discussions on Morse theory and polynomial optimization.  He would also like to  thank the anonymous referees for their   careful reading and detailed comments with many helpful suggestions.


\end{document}